From hfw@math.umn.edu Wed Apr  9 07:30:30 1997
Date: 8 Apr 1997 14:32:59 -0000
From: hfw@math.umn.edu
To: hans@math.wisc.edu
Subject: another revision

\magnification=\magstep1
\rightline{4-8-97 by HW}
\medskip
\centerline{{\bf On matrices for which norm bounds are attained.}}
\bigskip
\noindent Hans Schneider\footnote{$^*$}{Research supported by NSF Grant
DMS-9424346}
\smallskip                           
{\sl \noindent Department of Mathematics                      

\noindent Van Vleck Hall                             

\noindent 480 Lincoln Drive                        

\noindent University of Wisconsin-Madison       

\noindent Madison WI 53706 USA}
\smallskip
\noindent and
\smallskip
\noindent Hans F. Weinberger
\smallskip
{\sl \noindent School of Mathematics

\noindent University of Minnesota

\noindent 107 Vincent Hall

\noindent 206 Church Street

\noindent Minneapolis, MN, 55455 USA} 

\bigskip

\centerline{April 1997}

\vfil\eject

ABSTRACT
\medskip

Let $\|A\|_{p,q}$ be the norm induced on the matrix $A$ with $n$ rows
and $m$ columns by the H\"older $\ell_p$ and $\ell_q$
norms on $R^n$ and $R^m$
(or $C^n$ and $C^m$), respectively. It is easy to find an upper bound for
the ratio  $\|A\|_{r,s}/\|A\|_{p,q}$. In this paper we study the classes
of matrices for which the 
upper bound is attained. We shall show that for fixed $A$, attainment
of the bound 
depends only on the signs of $r-p$ and $s-q$. Various criteria depending on
these signs are obtained. For the special case $p=q=2$, the
set of all matrices for which the bound is attained is generated by
means of singular value decompositions.
\bigskip
\hrule  
\bigskip

\noindent 1. INTRODUCTION
\bigskip

Let $A$ be a matrix with $n$ rows and $m$
columns. 
If $A$ is considered as a complex transformation, let $\mu_1$ and
$\mu_2$ be norms on $C^m$, and let $\nu_1$ and
$\nu_2$ be norms on $C^n$. If $A$ is real and is considered as a
transformation from $R^m$ to $R^n$, let the $\mu_i$ be norms on $R^m$
and the $\nu_i$ be norms on $R^n$. Define the induced norms
$$
\|A\|^{(i)}=\max_{{\bf x}}\nu_i(A{\bf x})/\mu_i({\bf x})
$$ 
for $i=1$ and 2, where the maximum is taken over either $C^m$ or
$R^m$, as is appropriate. It was shown in [SS] (see also [HJ, p.303]) 
that 
$$
\|A\|^{(2)}\le \max_{{\bf x}}{\mu_1({\bf
x})\over\mu_2({\bf x})}\ \max_{{\bf y}}{\nu_2({\bf
y})\over\nu_1({\bf y})}\ \|A\|^{(1)},
\leqno(1.1) 
$$
and that equality is always attained for some $A\ne 0$. Here the
maxima are taken over $C^m$ and $C^n$ if $A$ is thought of as a
complex transformation, and over $R^m$ and $R^n$ if $A$ is a real and
its action is confined to $R^m$.

In this work we shall be concerned with characterizing the set of all
matrices $A$ for which equality is attained in (1.1), at least in some
cases. 

We shall show that this set can be described by the following
property.

{\bf Theorem 1.} {\sl If equality holds in the inequality (1.1),
then every maximizer ${\bf v}$ of the ratio $\nu_2(A{\bf
x})/\mu_2({\bf x})$ has the properties that 

\item{(i)} {\bf v} is also a maximizer
of the ratio $\mu_1({\bf
x})/\mu_2({\bf x})$, 

\item{(ii)} $A{\bf v}$ is a maximizer of the
ratio $\nu_2({\bf y})/\nu_1({\bf y})$, and

\item{(iii)} ${\bf v}$ is also a
maximizer of the ratio $\nu_1(A{\bf x})/\mu_1({\bf x})$.

Conversely, if there is one maximizer $\bf v$ of $\nu_1(A{\bf x})/\mu_1({\bf
x})$ which has the properties (i) and (ii), then equality holds in (1.1).}
\smallskip
Theorem 1 can only provide useful information if the two
maxima on the right and the corresponding maximizers are known.
Both of these conditions apply when the norms involved are H\"older
norms. We denote the $\ell_p$ norm by $\|\ \|_p$. For any $p$ and $q$
in the interval $[1,\infty]$ we define the induced norm
$$
\|A\|_{p,q}:=\max_{{\bf x}}{\|A{\bf x}\|_q\over\|{\bf x}\|_p}.
\leqno(1.2)
$$

The maximum of the ratio $\|{\bf x}\|_r/\|{\bf x}\|_p$ and the
corresponding maximizers are well known. In order to state the result
we recall that sgn$(z)$ is defined to be 1 if $z>0$, 0 if $z=0$, and
-1 if $z<0$, and that $[z]_+$ is defined to be $z$ if $z\ge0$ and 0 if
$z\le0$. We also define the three subsets of a real or complex vector
space of $m$-tuples or $n$-tuples.
$$
\eqalign{
&K_1=\{{\bf x}:\ {\rm all\ components\ of\ {\bf x}\ have\ equal\ absolute\
values}\}\cr
&K_{-1}=\{{\bf x}:\ {\rm at\ most\ one\ component\ of\ {\bf x}\ differs\
from\ 0}\}\cr
&K_0={\rm the\ whole\ vector\ space}.}
\leqno(1.3)
$$
The following result is found, e.g., in [HLP, p.~26 \#16 and
p.~29\#19]. 

{\bf Proposition 1.} For any ${\bf x}$ in $R^m$ or $C^m$ 
$$
\|{\bf x}\|_r\le m^{[(1/r)-(1/p)]_+}\ \|{\bf x}\|_p\qquad {\rm for}\ 
p,r\in[1,\infty].
\leqno(1.4)
$$
{\sl Equality holds if and only if the $m$-vector $\bf x$ lies in
$K_{{\rm sgn}(p-r)}$.}
\smallskip

By inserting Proposition 1 into the inequality (1.1) and into Theorem
1, we immediately obtain the following special case for the H\"older
spaces.

{\bf Proposition 2.} 
$$
\|A\|_{r,s}\le m^{[(1/p)-(1/r)]_+}\ n^{[(1/s)-(1/q)]_+}\ \|A\|_{p,q}
\qquad {\rm for}\ p,q,r,s \in [1,\infty].
\leqno(1.5)
$$

{\sl If equality holds in this inequality, then every maximizer ${\bf v}$ of
the ratio $\|A{\bf x}\|_s/\|{\bf x}\|_r$ has the properties 

\item{(i)} ${\bf
v}\in K_{{-\rm sgn}(p-r)}$, 

\item{(ii)} $A{\bf
v}\in K_{{\rm sgn}(q-s)}$, and

\item{(iii)} ${\bf v}$ is a maximizer of the
ratio $\|A{\bf x}\|_q/\|{\bf x}\|_p$.

Conversely, if there exists a maximizer ${\bf v}$ of the ratio
$\|A{\bf x}\|_q/\|{\bf x}\|_p$ which has the properties (i) and (ii),
then equality holds in (1.5).}
\smallskip

\noindent 
For the case $m=n$, $q=p$, $s=r$, the inequality (1.5)
was pointed out by Higham [H, p.124].  

\smallskip  \noindent
{\bf Remark.} {\sl 
The inequality (1.5) is equivalent to 
the following
monotonicity statement:} For fixed $s$,
$\|A\|_{r,s}$ is nondecreasing and $m^{1/r}\|A\|_{r,s}$ is
nonincreasing in $r$, and for fixed $r$, $n^{-1/s}\|A\|_{r,s}$
is nondecreasing and $\|A\|_{r,s}$ is nonincreasing in $s$.

{\sl The inequality (1.5) obviously implies the above statement. 
It is not difficult to show that the converse is true.
If, for instance, $p\le r$ and $q\le s$, then the
monotonicity statement implies that
$$
m^{1/r}\|A\|_{r,s}\le m^{1/p}\|A\|_{p,s}\le m^{1/p}\|A\|_{p,q},
$$
which implies the inequality (1.5) for this case.}
\smallskip

The trivial observation that for fixed $(p,q)$ the only dependence on
$(r,s)$ in Proposition 2 is through the functions sgn$(p-r)$ and
sgn$(q-s)$ immediately yields the following statement.

{\bf Proposition 3.} {\sl If equality holds in (1.5) , if
sgn$(p-r')$=sgn$(p-r)$, and if sgn$(q-s')$=sgn$(q-s)$, then equality
also holds in (1.5) when the pair $(r,s)$ is replaced by $(r',s')$.}
\smallskip

{\bf Remark.} {\sl By using the inequality (1.5) with $r=p'$ and
$s=q'$, one sees that Proposition 2 also shows that equality
in (1.5) implies that 
$$
\|A\|_{r',s'} =  m^{[(1/p')-(1/r')]_+}\ n^{[(1/s')-(1/q')]_+}\ \|A\|_{p',q'},
\qquad {\rm for}\ p,q,r,s \in [1,\infty],
$$
provided ${\rm  sgn}(p'-r')={\rm  sgn}(p-p')={\rm
sgn}(p-r)$ and 
${\rm  sgn}(q'-s')={\rm  sgn}(q-q')={\rm  sgn}(q-s)$.}
\smallskip

Proposition 3 shows that for a prescribed pair $(p,q)$ with
$1\le p,q\le\infty$, equality holds in (1.5) for some $r\ne p$ and $s\ne q$
if and only if $A$ lies in the appropriate one of at most four classes, which
we shall label by the extremal pair $(r,s)$, that is, the pair in
which each of these indices has the value 1 or $\infty$, for which
equality holds in (1.5). That is, we define
$$
\eqalign{
{\cal E}_{1,\infty}(p,q):=\{A:{\rm equality\ in\ (1.5)\ holds\ when}\  
r<p\ {\rm and}\ s>q\}, \cr
{\cal E}_{1,1}(p,q):=\{A:{\rm equality\ in\ (1.5)\ holds\ when}\  
r<p\ {\rm and}\ s<q\}, \cr
{\cal E}_{\infty,\infty}(p,q):=\{A:{\rm equality\ in\ (1.5)\ holds\ when}\  
r>p\ {\rm and}\ s>q\}, \cr
{\cal E}_{\infty,1}(p,q):=\{A:{\rm equality\ in\ (1.5)\ holds\ when}\  
r>p\ {\rm and}\ s<q\}.} 
\leqno(1.6)
$$
(Of course, when $p$ and or $q$ has one of the extreme values, some of
these classes are trivial.)  This work is concerned with
characterizing the members of these four classes.

When $p=q=2$, Proposition 2 enables us to give a characterization of
all matrices for which equality holds in the bound (1.5). As usual, we
denote the Hermitian transpose of a matrix $A$ by $A^*$.

{\bf Theorem 2.} {\sl If $r,s\in[1,\infty]$, the equality
$$
\|A\|_{r,s}= m^{[(1/2)-(1/r)]_+}\ n^{[(1/s)-(1/2)]_+}\ \|A\|_{2,2}.
$$
is valid if and only if $A$ has a singular value decomposition
$$
A=U\Sigma V^*
$$
in which 

\item{(i)} the first column of the unitary matrix $U$ is in $K_{{\rm
sgn}(2-s)}$, 

\item{(ii)} the first column of the unitary matrix $V$ is in
$K_{-{\rm sgn}(2-r)}$, and

\item{(iii)} the (11) entry of the nonnegative
diagonal matrix $\Sigma$ is its maximal entry.}
\smallskip

The first two theorems will be proved in Section 2.

When $p$ and $q$ are not both 2, Proposition 2 will will be used to
obtain characterizations of the classes in (1.6).  Our most complete
characterization is for the class ${\cal E}_{1,\infty}(p,q)$, which is
treated in Section 3.

{\bf Theorem 3.} {\sl Let $\rho$ denote the largest absolute value of
the entries of $A$. 
 
Then $\|A\|_{1,\infty}=\rho$, and $A\in{\cal E}_{1,\infty}(p,q)$ if
and only if $A$ has the properties

\item{(i)} every
entry of $A$ which has the absolute value $\rho$ is the only nonzero
element of its row and of its column, and

\item{(ii)} if $C$ is the matrix
obtained from $A$ by replacing all elements of absolute value $\rho$
by zero, then $\|C\|_{p,q}\le\rho$.

If $p>q$, then 
$A\in{\cal E}_{1,\infty}(p,q)$ 
if and only if $A$ has at most one nonzero entry.}
\smallskip

\noindent Theorem $3'$ in Section 3 shows that the Property (i) is sufficient
for the existence of a $p>1$ and a $q<\infty$ such that $A\in
{\cal E}_{1,\infty}(p,q)$.

Section 4 deals with the cases in which $r<p$ and $s<q$ or $r>p$ and
$s>q$. We shall establish the following results.

{\bf Theorem 4.} {\sl Let $\sigma$ denote the largest $\ell_1$ norm of
the columns of $A$, so that $\sigma=\|A\|_{1,1}$.  

If $A\in{\cal E}_{1,1}(p,q)$ then $A$ has the properties

\item{(i)} the
entries of any column whose $\ell_1$ norm is equal to $\sigma$ all
have the same absolute value $n^{-1}\sigma$, 

\item{(ii)} every column with
this property is orthogonal to all the other columns of $A$, and

\item{(iii)} $\sigma=n^{1-(1/q)}\|A\|_{p,q}$.

Conversely, if the matrix $A$ has a column all of whose entries have
the absolute values $n^{-1/q}\|A\|_{p,q}$, then $A\in{\cal
E}_{1,1}(p,q)$. 

If $p>2$, then $A\in{\cal E}_{1,1}(p,q)$ if and only if $A$ has only
one nonzero column, and all the entries of this column have the same
absolute value.}
\smallskip

{\bf Theorem 5.} {\sl Let $\sigma$ denote the largest $\ell_1$ norm of
the rows of $A$, so that $\sigma=\|A\|_{\infty,\infty}$.  

If $A\in{\cal E}_{\infty,\infty}(p,q)$, then $A$ has the properties

\item{(i)} the entries of any row whose $\ell_1$ norm is equal to
$\sigma$ all have the same absolute value $m^{-1}\sigma$, 

\item{(ii)} every
row with this property is orthogonal to all the other rows of $A$, and

\item{(iii)} $\sigma=m^{1/p}\|A\|_{p,q}$.

Conversely, if the matrix $A$ has a row all of whose entries have
the absolute values $m^{(1/p)-1}\|A\|_{p,q}$, then $A\in{\cal
E}_{\infty,\infty}(p,q)$. 

If $q<2$, then $A\in{\cal E}_{\infty,\infty}(p,q)$ if and only if $A$
has only one nonzero row, and all the entries of this row have the
same absolute value.}
\smallskip

Theorem $4'$ in Section 4 shows that the Properties (i) and (ii) of
Theorem 4 are sufficient for the existence of $p>1$ and $q>1$ such
that $A\in{\cal E}_{1,1}(p,q)$. Analogously, Theorem 5' states that
the properties 
(i) and (ii) of Theorem 5 imply the existence of finite $p$ and $q$
such that $A\in{\cal E}_{\infty,\infty}(p,q)$

Section 5 considers the case where $r>p$ and $s<q$. The following
result is obtained.

{\bf Theorem 6.} {\sl $A\in{\cal E}_{\infty,1}(p,q)$ if and only if
there is a vector $\bf v$ with the properties

\item{(i)} ${\bf v}$ is an
eigenvector of the matrix $A^*A$, 

\item{(ii)} all the entries of $\bf v$ have
the absolute value 1, 

\item{(iii)} all the entries of $A{\bf v}$ have the
same absolute value $\tau$, and 

\item{(iv)} $\tau=m^{1/p}n^{-1/q}\|A\|_{p,q}$.

In particular, $A\in{\cal E}_{\infty,1}(2,2)$
if and only if $A^*A$ has an eigenvector $\bf v$ with the properties (ii) and
(iii) which corresponds to its largest eigenvalue.}
\smallskip

We observe that when the matrix $A$ is real, one has a choice of
defining the induced norm $\|A\|_{p,q}$ with respect to either the
real or the complex H\"older spaces, and that these two norms may differ
for some $(p,q)$. Our results are valid for either choice.

Consider, for instance, the matrix $A=\pmatrix{1&1\cr-1&1}$. The last
statement of Theorem 6 with the complex eigenvector $(1,\ i)$ of
$A^*A=2I$, shows that when $r\ge2\ge s$ the norm $\|A\|_{r,s}$ on the
complex vector space $C^2$ is equal to $2^{(1/s)-(1/r)+(1/2)}$. On the
other hand, a simple computation shows that on the real vector spaces,
$\|A\|_{\infty,1}=2$ while $\|A\|_{2,2}$ is still $2^{1/2}$. Thus in
the real norm, equality does not hold in (1.5) when $p=q=2$, $r>2$, and
$s<2$. Therefore the
real norm $\|A\|_{r,s}$ is strictly less than $2^{(1/s)-(1/r)+(1/2)}$,
and hence less than the complex norm, when $r>2$ and $s<2$.
\bigskip

\noindent 2. PROOFS OF THEOREMS 1 AND 2. 
\bigskip
 
We begin by proving Theorem 1.

{\sl Proof of Theorem 1.} We recall the derivation in [SS] of the
inequality (1.1). For any $\bf x\ne 0$ with $A{\bf x}\ne 0$ we have
$$
{\nu_2(A{\bf x})\over\mu_2({\bf x})}={\mu_1({\bf x})\over\mu_2({\bf
x})}{\nu_2(A{\bf x})\over\nu_1(A{\bf
x})}{\nu_1(A{\bf x})\over\mu_1({\bf x})}.
\leqno(2.1)
$$
Because the maximum of a product of nonnegative numbers is bounded by
the product of the maxima, we obtain the inequality (1.1).

Suppose there is a maximizer of the left-hand side of
(2.1), which is not a maximizer of one of the factors on the right.
Since all the factors are bounded by their maxima and one of them is
strictly less than its maximum,
the right-hand side of (1.1) is strictly greater than the left-hand side. Therefore 
the condition of Proposition 1 is necessary for equality.

If there is a maximizer $\bf v$ of all three quotients on the right of
(2.1), then the maximum of the left-hand side is bounded below by the
right-hand side of (1.1). Since we already know that it is bounded above
by the same quantity, we conclude that equality holds in (1.1). This
establishes Theorem 1.
\smallskip

{\sl Proof of Theorem 2.} We observe that a maximizer of the ratio
$\|A{\bf v}\|_2/\|{\bf v}\|_2$ is an eigenvector of the matrix $A^*A$
which corresponds to its largest eigenvalue. By Proposition 2,
equality in (1.5) with $p=q=2$ implies that a
maximizer $\bf v$ of $\|A{\bf x}\|_r/\|{\bf x}\|_q$ is such an
eigenvector, that it is in $K_{-{\rm sgn}(2-r)}$, and that the eigenvector
$A{\bf v}$ of $AA^*$ is in $K_{{\rm sgn}(2-s)}$.

Thus we can construct (see. e.g., the proof of Theorem 2.3-1 in [GVL])
a singular value decomposition $A=U\Sigma V^*$ in which the first
column of the unitary matrix $U$ is the vector $\|A{\bf
v}\|_2^{-1}A{\bf v}\in K_{{\rm sgn}(2-s)}$ and the first row of the
unitary matrix $V$ is $\|{\bf v}\|_2^{-1}{\bf v}\in K_{-{\rm
sgn}(2-r)}$. The (11) element of the nonnegative diagonal matrix
$\Sigma$ is the square root of the largest eigenvalue of $A^*A$, which
is the maximal element of $\Sigma$.

The converse follows from the fact that the first column of $V$ is a
maximizer of $\|A{\bf x}\|_2/\|{\bf x}\|_2$ and the converse statement
of Proposition 2, so that Theorem 2 is proved. 

{\bf Remark.} {\sl If the matrix $A$ is a scalar multiple of a unitary
matrix and the absolute values of all its entries are equal to a
number $\rho$, then $A$ has a singular value decomposition with
$U=n^{-1/2}\rho^{-1}A$, $\Sigma=n^{1/2}\rho I$, and $V=I$, and another
singular value decomposition with $U=I$, $\Sigma=n^{1/2}\rho I$, and
$V=n^{-1/2}\rho^{-1}A^*$. Hence Theorem 2 shows that 
$A$ lies in both ${\cal E}_{1,1}(2,2)$ and ${\cal E}_{\infty,\infty}(2,2)$.

Examples of such matrices include the Hadamard matrices, which are
orthogonal matrices whose entries have the values $\pm1$ (see [H,
p.~128, \S6.13]), and the matrices which represent the finite
Fourier transforms.}

\bigskip
\noindent 3. THE CLASS ${\cal E}_{1,\infty}(p,q)$.
\bigskip

 The following lemma will be used in the
proofs of Theorems 3, 4, and 6. We recall the definition of the
conjugate index $p^*=p/(p-1)$ of an index $p$, and the fact that $A^*$
denotes the Hermitian transpose of the matrix $A$.

We also recall the identity
$$
\|A^*\|_{q^*,p^*}=\|A\|_{p,q},
\leqno(3.1)
$$
which simply states that the norm of the adjoint of a transformation
is equal to the norm of the transformation.

{\bf Lemma 3.1.} {\sl Suppose that a maximizer $\bf v$ of $\|A{\bf
x}\|_q/\|{\bf x}\|_p$ has the properties that 

\item{(i)} all its nonzero
components have the same absolute value, and 

\item{(ii)} the same is true of
$A{\bf v}$. 

\noindent If $1<p<\infty$, or $p=1$ and ${\bf v}\in K_1$, or
$p=\infty$ and ${\bf v}\in K_{-1}$, then $\bf v$ is an eigenvector of
the matrix $A^*A$.}

{\sl Proof.} Because of the duality relation (3.1), we have
$$
\eqalign{
\|A^*A{\bf v}\|_{p^*}&\le\|A\|_{p,q}\|A{\bf v}\|_{q^*} \cr
  &=(\|A{\bf v}\|_{q}/\|{\bf v}\|_p)\|A{\bf v}\|_{q^*}.}
$$
It is easily seen from the property (ii) that
$$
\|A{\bf v}\|_2^2=\|A{\bf v}\|_{q}\|A{\bf v}\|_{q^*}.
$$
Therefore
$$
{\bf v}\cdot A^*A{\bf v}=\|A{\bf v}\|_{q}\|A{\bf v}\|_{q^*}\ge\|{\bf
v}\|_p\|A^*A{\bf v}\|_{p^*}.
$$
This shows that equality holds in the H\"older inequality for the
bilinear form ${\bf v}\cdot A^*A{\bf v}$ in
$\ell_p\times\ell_{p^*}$. If $1<p<\infty$, this implies that the
vector $A^*A{\bf v}$ must be a multiple of the vector with components
$\|{\bf v}\|_p^{p-2}v_j$. (See, e.g., [HLP p.~26\#14].) By Property (i) this vector is a multiple
of $\bf v$, which proves the result for this case.

If $p=1$ so that $p^*=\infty$, and if $\bf v$ has no zero component,
it is easily seen that equality in H\"older's inequality implies that
$A^*A{\bf v}$ is proportional to the vector with components
$|v_j|^{-1}v_j$, and we reach the same conclusion. This is the case when
$p=1$ and ${\bf v}\in K_1$. 

If $p=\infty$ so that $p^*=1$, one easily sees that equality in the
H\"older inequality implies that $A^*A{\bf v}$ has zero components
where $\bf v$ does. Therefore,if ${\bf v}\in K_{-1}$ so that it has
only one nonzero component, $A^*A{\bf v}$ is again proportional to
$\bf v$.

Thus the Lemma is proved in all cases.
\smallskip

{\sl Proof of Theorem 3.} It is easily
verified that $\|A{\bf x}\|_\infty/\|{\bf x}\|_1\le \rho$, the largest
absolute value of any entry of $A$, and that this bound is attained
when $\bf x$ is in the direction of a coordinate which corresponds to
a column in which an element of magnitude $\rho$ occurs. Thus
$\|A\|_{1,\infty}=\rho$, which is the first statement of the Theorem.

Moreover, a unit coordinate 
vector $\bf v$ in the direction of a column which contains an element
of magnitude $\rho$ is a maximizer of the ratio.

Suppose now that $A\in{\cal E}_{1,\infty}(p,q)$.  Then equality in
(1.5) holds for $r=1$ and $s=\infty$.  Proposition 2 shows that if
$\bf v$ is a unit vector in the direction of a column of $A$ with a
maximal element, this column has exactly one nonzero element, and $\bf
v$ is a maximizer of the ratio $\|A{\bf x}\|_q/\|{\bf x}\|_p$.  The
first of these properties says that any column of $A$ which contains
an element of magnitude $\rho$ has but one nonzero element, while the
second property implies that the absolute value $\rho$ of the nonzero
element equals $\|A\|_{p,q}=\|A\|_{1,\infty}$. There may, of course,
be several maximizers, and therefore several columns with singleton
elements of magnitude $\rho$.

Since $\bf v$ and $A{\bf v}$ are both in coordinate directions and
$p>1$, we can apply Lemma 3.2 to show that $\bf v$ is an eigenvector
of $A^*A$. Therefore if $\bf x$ is a coordinate vector orthogonal to
$\bf v$, it is also 
orthogonal to $A^*A{\bf v}$, which implies that $A{\bf x}$ is
orthogonal to $A{\bf v}$. This means that a column which contains a
single nonzero element of magnitude
$\rho$ is orthogonal to all the other columns of $A$. In other words,
an  element of magnitude $\rho$ is also the only nonzero element of
its row as well as of its column, so that Property (i) is established.

If we choose a trial vector $\bf x$ whose components in the directions
of the columns with elements of magnitude $\rho$ are zero, then
$A{\bf x}=C{\bf x}$ where $C$ is defined in the statement of Theorem
1. Therefore $\|C\|_{p,q}\le\|A\|_{p,q}=\rho$. This is Property (ii).

To prove the converse statement for $p\le q$, we define $B=A-C$, and
decompose any vector $\bf x$ into ${\bf y}+{\bf z}$, where the
components of $\bf z$ are zero in the directions corresponding to
columns which contain elements of magnitude $\rho$ and the components
of $\bf y$ in the remaining directions vanish. Then by Property (i)
and two applications of Proposition 1
$$
\eqalign{
\|A{\bf x}\|_q&=\{(\rho\|{\bf y}\|_q)^q+\|C{\bf z}\|_q^q\}^{1/q}\cr
    &\le\{(\rho\|{\bf y}\|_p)^q+(\|C\|_{p,q}\|{\bf z}\|_p)^q\}^{1/q}\cr
    &\le\{(\rho\|{\bf y}\|_p)^p+(\|C\|_{p,q}\|{\bf z}\|_p)^p\}^{1/p}\cr
    &\le\max\{\rho,\|C\|_{p,q}\}\|{\bf x}\|_p.}
$$
That is,
$$
\|A\|_{p,q}=\max\{\rho,\|C\|_{p,q}\}.
\leqno(3.2)
$$
Thus Property (ii) shows that $\|A\|_{p,q}=\rho=\|A\|_{1,\infty}$,
and the proof of the converse statement is complete.

To prove the last assertion of Theorem 3
assume that $p>q$ and that $A\in{\cal E}_{1,\infty}(p,q)$.
Choose a trial vector $\bf x$ whose component in
the direction of a column with a singleton element of magnitude
$\rho$ is one and which has one other nonzero component $\alpha$. Let $b$
be any entry of $A$ in the column which corresponds to $\alpha$. Then
because $p>q$,
$$
{\|A{\bf x}\|_q\over\|{\bf x}\|_p}
\ge {(\rho^q+|\alpha b|^q)^{1/q}\over(1+|\alpha|^p)^{1/p}}
=\rho+(1/p)\rho^{1-q}|b|^q|\alpha|^q+o(|\alpha|^q)
$$
for small $\alpha$. Because $\rho=\|A\|_{p,q}$, the right-hand side
must be bounded by $\rho$, and we conclude that $b=0$. Because $b$ is
an arbitrary element of any column other than that with the entry of
magnitude $\rho$, we conclude that all other columns of $A$ are
zero, so that $A$ has only one nonzero entry.

Finally, a simple computation shows that if $A$ has only one nonzero
entry, and if the magnitude of this entry is $\rho$, then
$\|A\|_{r,s}=\rho$ for all $r$ and $s$, so that $A\in{\cal
E}_{1,\infty}(p,q)$. 

Thus all parts of Theorem 3 have been established.
\smallskip 

Because it is difficult to compute the $p,q$ norm for most $p$ and
$q$, it is difficult to verify Property (ii) of Theorem 3. We shall
show that the easily verified Property (i) is sufficient to assure the
existence of some $p>1$ and $s<\infty$ such that equality holds in
(1.5) when $r<p$ and $s>q$.

{\bf Theorem $3'$.} {\sl Let $A$ have Property (i) of Theorem 3.
Let $C$ be the matrix obtained from $A$ by replacing all elements
of absolute value $\rho=\|A\|_{1,\infty}$ by 0, so that
$\|C\|_{1,\infty}<\rho$. If $p$ and $q$ satisfy the inequalities
$p\le q$ and
$$
m^{1-(1/p)}n^{(1/q)}\|C\|_{1,\infty}\le \rho,
\leqno(3.3)
$$
then $A\in{\cal E}_{1,\infty}(p,q)$. The inequality (3.3) is satisfied
if $p$ is sufficiently close to $1$ and $q$ is sufficiently large.}

\smallskip
  
{\sl Proof.} Since (1.5) shows that 
$$
\|C\|_{p,q}\le
m^{1-(1/p)}n^{1/p}\|C\|_{1,\infty},
$$
the inequality (3.3) and the equation (3.2) imply that
$\|A\|_{p,q}=\rho$. That is, Property (ii) of Theorem 3 holds, and
the conclusion $A\in{\cal E}_{1,\infty}(p,q)$ follows.
\bigskip

\noindent 4. THE CLASSES ${\cal E}_{1,1}(p,q)$ AND 
${\cal E}_{\infty,\infty}(p,q)$.
\bigskip

{\sl Proof of Theorem 4.} Suppose that $A\in{\cal E}_{1,1}(p,q)$, so
that equality holds
in (1.5) with $r=s=1$.
The triangle inequality shows that $\|A{\bf x}\|_1/\|{\bf x}\|_1\le
\sigma$, the largest $\ell_1$ norm of the columns of $A$. Moreover,
this bound is attained when $\bf x$ is in the direction of any
coordinate whose corresponding column has the $\ell_1$ norm
$\sigma$. Thus if $\bf v$ is a coordinate vector in such a direction,
it is a maximizer for the ratio.

Proposition 2 states that if
$\bf v$ is a unit vector in one of these coordinate directions, the
elements of the corresponding column $A{\bf v}$ must have equal
absolute values, and $\bf v$ must also be a maximizer of $\|A{\bf
x}\|_q/\|{\bf x}\|_p$. These two facts give the properties (i) and
(iii) of Theorem 3.

Since $p>1$, Lemma 3.1 shows that $\bf v$ is an eigenvector of
$A^*A$. As in the proof of Theorem 3, this implies that if $\bf x$ is
a coordinate vector perpendicular to $\bf v$, then it is also
perpendicular to $A^*A{\bf v}$, so that the column $A{\bf x}$ is
perpendicular to $A{\bf v}$. This is the property (ii)

To prove the converse statement, we observe that if $A$ has a column
whose elements have the absolute value $n^{-1/q}\|A\|_{p,q}$, then
a unit vector $\bf v$ in the direction of this column is a maximizer
of the ratio $\|A{\bf x}\|_q/\|{\bf x}\|_p$. Therefore the converse
statement of Proposition 2 implies that equality holds in (1.5), and
hence that $A\in{\cal E}_{1,1}(p,q)$.

To prove the last statement of Theorem 4, we suppose that $A\in{\cal
E}_{1,1}(p,q)$, and that $p>2$.
Then there is at least one column $\bf c$ of $A$ all
of whose entries have the absolute value
$n^{-1}\sigma=n^{-1/q}\|A\|_{p,q}$.  Let $\bf c$ be one such column,
let $\bf b$ be any other column of $A$, and let $\alpha$ be a real
parameter. The adjoint relation (3.1) leads to the inequality
$$
\|A^*({\bf c}+\alpha{\bf b})\|_{p^*}\le\|A\|_{p,q}\|{\bf
c}+\alpha{\bf b}\|_{q^*}=n^{-1+(1/q)}\sigma\|{\bf
c}+\alpha{\bf b}\|_{q^*}.
\leqno(4.1)
$$

We observe
that for small $\alpha$
$$
|c_j+\alpha b_j|^{q^*}=|c_j|^{q^*}+q^*\ {\rm
    Re}(\alpha|c_j|^{q^*-2}\bar{c_j}b_j)+O(\alpha^2).
$$
We sum on $j$ and use the properties that
the entries of $\bf c$ all have the absolute value
$n^{-1}\sigma$ and that $\bf b$ is orthogonal to $\bf c$ to find
that 
$$
\|{\bf c}+\alpha{\bf
b}\|_{q^*}^{q^*}=n^{-1+(1/q^*)}\sigma+O(\alpha^2)=n^{-1/q}\sigma+O(\alpha^2).
\leqno(4.2)
$$ 

Since $\bf c$ and $\bf b$ are orthogonal, the entry of $A^*({\bf
c}+\alpha{\bf b})$ which corresponds to the column $\bf c$ is
$n^{-1}\sigma^2$, while the entry which corresponds to
the column $\bf b$ is $\alpha\|{\bf b}\|_2^2$. We 
obtain a lower bound for the
left-hand side of (4.1) by replacing all the other entries by
zero. For small $\alpha$ this lower bound takes the form
$$
\|A^*({\bf c}+\alpha{\bf b})\|_{p^*}
\ge n^{-1}\sigma^2+(p^*)^{-1}(n^{-1}\sigma^2)^{1-p^*}\|{\bf
b}\|_2^{2p^*}\alpha^{p^*}
+O(\alpha^{2p^*}).
$$
where $K>0$.

By putting this and (4.2) into (4.1), we find the inequality
$$
n^{-1}\sigma^2+(p^*)^{-1}(n^{-1}\sigma^2)^{1-p^*}\|{\bf
b}\|_2^{2p^*}\alpha^{p^*}
+O(\alpha^{2p^*})\le n^{-1}\sigma^2+O(\alpha^2).
$$
We observe that $p^*<2$ because $p>2$. We cancel the first terms from
the two sides, divide by $\alpha^{p^*}$, and let $\alpha$ approach
zero to see that $\|{\bf b}\|_2=0$. That is, every column other than
$\bf c$ is zero. This establishes the last statement of Theorem 4, and
the Theorem is proved.
\smallskip

Theorem 5 will follow easily from Theorem 4 and the following lemma.

{\bf Lemma 4.1.}  {\sl If equality holds in (1.5), then equality also
holds when $A$ is replaced by $A^*$, the pair $(r,s)$ is replaced by
$(s^*,r^*)$, and the pair $(p,q)$ is replaced by $(q^*,p^*)$.}

{\sl Proof.} We recall the adjoint equation (3.1), namely 
$\|A^*\|_{q^*,p^*}=\|A\|_{p,q}$ 
We also note that in going
from $A$ to $A^*$ the dimensions $m$ and $n$ are interchanged, and that
by definition $(1/q^*)-(1/s^*)=(1/s)-(1/q)$ and
$(1/r^*)-(1/p^*)=(1/p)-(1/r)$. Therefore, the replacements indicated
in the Lemma leave both sides of (1.5) unchanged, which proves the
Lemma.
\smallskip

{\sl Proof of Theorem 5.} By Lemma 4.1, $A\in{\cal
E}_{\infty,\infty}(p,q)$ if and only if $A^*\in{\cal
E}_{1,1}(q^*,p^*)$.  Since $s>q$ implies $s^*<q^*$ and $r>p$ implies
$r^*<p^*$, the application of Theorem 4 to $A^*$ with the above index
replacements gives the statement of Theorem 5.
\smallskip

As in the case of Theorem 3, it is difficult to verify the last
hypothesis of Theorem 4. We shall prove that the easily
verified Properties (i) and (ii) are sufficient to assure the
existence of $p,q\in(1,\infty]$ such that 
$A\in{\cal E}_{1,1}(p,q)$

{\bf Theorem $4'$.} {\sl Let $A$ have the Properties (i) and (ii) of
Theorem 4. Let $C$ be the matrix obtained from $A$ by replacing all
columns with the $\ell_1$ norm $\sigma=\|A\|_{1,1}$ by zero, so that
$\|C\|_{1,1}<\sigma$. If $p\le2$ satisfies the inequality 
$$
\eqalign{
(2mn&)^{1-(1/p)}\sigma^{-1}\|C\|_1^0 
  +[2^{1-(1/p)}-1]\cr
&\cdot[(p/2)^{1/(2-p)}n^{(-3p^2+2p+4)/[2p(2-p)]}m^{-2(p-1)/(2-p)}
(\|C\|_{1,1}/\sigma)^{2/(2-p)}]\le 1,}
\leqno(4.3)
$$
and $q\le p$, then $A\in{\cal E}_{1,1}(p,q)$.  The inequality (4.3) is
satisfied when $p$ is sufficiently near 1.}
\smallskip

{\sl Proof of Theorem $4'$.} We recall that $C$ is the matrix obtained
from $A$ by replacing those columns whose $\ell_1$ norm is $\sigma$ by
0. Thus $\|C\|_{1,1}<\sigma$. Let $B=A-C$, so that all the nonzero
elements of $B$ have the magnitude $n^{-1}\sigma$, and every column of
$B$ is orthogonal to all other columns of $A$. To establish the
Theorem, we only need to show that the inequality (4.3) implies that
$\|A\|_{p,p}=n^{-1+(1/p)}\sigma=n^{-1+(1/p)}\|A\|_{1,1}$.

Decompose an arbitrary vector ${\bf x}\ne0$ into ${\bf x}={\bf y}+{\bf
z}$, where ${\bf z}$ is obtained from ${\bf x}$ by replacing those
elements which correspond to the nonzero columns of $B$ by zero, and
${\bf y}={\bf x}-{\bf z}$.
 
We see from the conditions (i) and (ii) of Theorem 4 that for the
above decomposition ${\bf x}={\bf y}+{\bf z}$,
$$
\|B{\bf x}\|_2^2=n^{-1}\sigma^2\|{\bf y}\|_2^2.
\leqno(4.4)
$$
In particular, $\|B\|_{2,2}=n^{-1/2}\sigma$, so that $B$ satisfies the
conditions of Theorem 3 with $p=q=2$. Therefore,
$$
\|B\|_{r,s}=n^{-1+(1/s)}\sigma
\leqno(4.5)
$$
for all $r$ and $s$ in the interval [1,2].  On the other hand, the
inequality (1.5) shows that
$$
\|C\|_{r,s}\le m^{1-(1/r)}\|C\|_{1,1}.
\leqno(4.6)
$$
for $r, s\ge1$.
Therefore 
if $p\le2$, the triangle inequality shows that
$$
\|A{\bf x}\|_p\le n^{-1+(1/p)}\sigma\|{\bf y}\|_p+m^{1-(1/p)}\|C\|_{1,1}
\|{\bf z}\|_p.
\leqno(4.7)
$$

Proposition 2 shows that 
$$
\|{\bf x}\|_p=(\|{\bf y}\|_p^p+\|{\bf z}\|_p^p)^{1/p}
  \ge 2^{-1+(1/p)}(\|{\bf y}\|_p+\|{\bf z}\|_p).
\leqno(4.8)
$$

We see from (4.7) and (4.8) that
$$
{\|A{\bf x}\|_p\over \|{\bf x}\|_p}\le n^{-1+(1/p)}\sigma
\leqno(4.9)
$$
whenever
$$ 
\|{\bf y}\|_p\le{1-(2mn)^{1-(1/p)}\sigma^{-1}\|C\|_1^0\over
2^{1-(1/p)}-1}\|{\bf z}\|_p.
\leqno(4.10)
$$

Thus the bound (4.4) is valid when the ratio
$\|{\bf y}\|_p/\|{\bf z}\|_p$ is not too large. To obtain this bound
for larger values of this ratio, we note that
$$
\|A{\bf x}\|_p^p=\sum_{j=1}^n(|(B{\bf y})_j|^2+2{\rm
Re}[\overline {(B{\bf y})_j} 
(C{\bf z})_j]+|(C{\bf z})_j|^2)^{p/2}.
\leqno(4.11)
$$

Because $p\le2$, the function $w^{p/2}$ is concave, so that for any
positive $d$ and $w$
$$
w^{p/2}\le
d^{p/2}+(p/2)d^{(p/2)-1} (w-d).
$$
We apply this inequality with $d=n^{-1}\|B{\bf y}\|_2^2$ to each term of
the sum on the right of (4.11) and use the fact that the range of $C$
is orthogonal to the range of $B$ to see that
$$
\|A{\bf x}\|_p^p\le n^{1-(p/2)}\|B{\bf y}\|_2^p
   +(p/2)n^{1-(p/2)}\|B{\bf y}\|_2^{p-2}\|C{\bf z}\|_2^2.
\leqno(4.12)
$$

The equation (4.5) shows that the first term on the right is bounded
by \break $n^{-1+(1/p)}(\sigma\|{\bf y}\|_p)^p$. Therefore we see that
the inequality (4.9) is valid when
$$
(p/2)n^{1-(p/2)}\|B{\bf y}\|_2^{p-2}\|C{\bf z}\|_2^2\le n^{p-1}\sigma^p
 \|z\|_p^p.
\leqno(4.13)
$$
We see from (4.6) that $\|C{\bf z}\|_2\le m^{1-(1/p)}\|C\|_{1,1}\|{\bf z}\|_p$, and 
from (4.4) and (1.4) that 
$$
\|B{\bf y}\|_2=n^{-1/2}\sigma\|{\bf y}\|_2\ge n^{-1/p}\sigma\|{\bf y}\|_p.
$$
Therefore the inequality (4.13), and hence also (4.9), is implied by 
$$
\|{\bf y}\|_p\ge
(p/2)^{1/(2-p)}n^{(-3p^2+2p+4)/[2p(2-p)]}m^{-2(p-1)/(2-p)}
(\|C\|_{1,1}/\sigma)^{2/(2-p)}\|{\bf z}\|_p.
\leqno(4.14)
$$

We now observe that the inequality (4.3) states that the coefficient
on the right of (4.14) is no larger than that in (4.10). Therefore at
least one of these inequalities inequalities is satisfied for every
$\bf y$ and $\bf z$. That is, the inequality (4.9) holds for all $\bf
x$, so that $\|A\|_{p,p}\le
n^{-1+1/p}\sigma=n^{-1+1/p}\|A\|_{1,1}$. Because (1.5) gives the
inequality in the opposite direction, we conclude that 
$A\in{\cal E}_{1,1}(p,p)$, and hence also that $A\in{\cal E}_{1,1}(p,q)$
for any $q\le p$. Thus Theorem $4'$ is established.
\smallskip 

By using Lemma 4.1 and applying Theorem $4'$ to $A^*$, we obtain the
analogous result.

{\bf Theorem $5'$.} {\sl Let $A$ have the Properties (i) and (ii) of
Theorem 5. Let $C$ be the matrix obtained from $A$ by replacing all
rows with the $\ell_1$ norm $\sigma=\|A\|_{1,1}$ by zero, so that
$\|C\|_{1,1}<\sigma$. If $q\ge2$ satisfies the inequality 
$$
\eqalign{
(2mn&)^{1/q}\sigma^{-1}\|C\|_1 
  +[2^{1/q}-1]\cr
&\cdot[(q^*/2)^{1/(2-q^*)}
m^{(-3(q^*)^2+2q^*+4)/[2q^*(2-q^*)]}n^{-2(q^*-1)/(2q^*p)}
(\|C\|_{1,1}/\sigma)^{2/(2-q^*)}]\le 1,}
\leqno(4.15)
$$
and $p\ge q$, then $A\in{\cal E}_{\infty,\infty}(p,q)$. The inequality
(4.15) is satisfied when 
$q$ is sufficiently large.}

\bigskip
\noindent 5. THE CLASS ${\cal E}_{\infty,1}(p,q)$.
\bigskip

{\sl Proof of Theorem 6.} Suppose that $A\in{\cal E}_{\infty,1}(p,q)$.
Proposition 2 shows that every maximizing vector
${\bf v}$ of the ratio $\|A{\bf x}\|_1/\|{\bf x}\|_\infty$ has the
properties (ii) its components have equal absolute values, which we
normalize to 1; (iii) the components of $A{\bf v}$ have equal absolute
values, which we call $\tau$; and (iv) $\|A{\bf v}\|_q/\|{\bf
v}\|_p=\|A\|_{p,q}$. Because $p<\infty$, Lemma 3.1 shows that $\bf v$
is an eigenvector of $A^*A$, which is Property (i). Thus the first
part of Theorem 6 is proved.

On the other hand, a vector $\bf v$ with the properties (ii), (iii),
and (iv) is a maximizer of the ratio $\|A{\bf x}\|_q/\|{\bf x}\|_p$,
so that Proposition 2 also establishes the converse statement.

The last statement of  
Theorem 6 clearly follows from the rest when $p=q=2$, so the Theorem
is proved. 

\smallskip

We are unable to find an analog of Theorems $3'$, $4'$, and $5'$ for this
case. We confine ourselves to the following simple observations.

1. If we define $V$ to be the diagonal unitary matrix whose diagonal
entries are the components of $\bf v$ and $D$ to be the diagonal
unitary matrix whose diagonal entries are the components of the vector
$\tau^{-1}\overline{A{\bf v}}$, the conditions of Theorem 4 imply that
all the row sums of the matrix $DAV$ are $\tau$ and that all its
column sums are $m^{-1}n\tau$. Conversely, if one can find two
matrices $D$ and $V$ with these properties, then the vector $\bf v$
whose components are the diagonal entries of $V$ has the properties
(i), (ii), and (iii) of Theorem 4. Thus equality holds in (1.5) for $r
> p$ and $s<q$ if and only if there are matrices $D$ and $V$ with
these properties and $\tau=m^{1/p}n^{-1/q}\|A\|_{p,q}$.

2. A sufficient condition for $A\in{\cal E}_{\infty,1}(2,2)$
is that there exist diagonal unitary matrices $D$
and $V$ such that the matrix $DAV$ has nonnegative entries, equal row
sums, and equal column sums. When $m=n$, $DAV$ is a multiple of
a doubly stochastic matrix.

3. The matrices with a single nonzero element which occur in the last
statement of Theorem 3 can be thought of as the tensor product of two
vectors in $K_{-1}$. Similarly, the matrices in the last statements of
Theorems 4 and 5 are tensor products. It is easily verified that if
$A={\bf c}\otimes{\bf b}$ so that its entries have the form $c_ib_j$,
then $\|A\|_{r,s}=\|{\bf b}\|_{r^*}\|{\bf c}\|_s$. Then Proposition 1
shows that when $A={\bf c}\otimes{\bf b}$, equality holds in (1.5)
if and only if ${\bf a}\in K_{-{\rm sgn}(p-r)}$ and ${\bf c}\in
K_{{\rm sgn}(q-s)}$. 

Theorem 6 and the fact that $\|A\|_{1,\infty}=\rho$ show that
$A\in{\cal E}_{\infty,1}(p,q)$ for all $p\in[1,\infty)$ and $q\in(1,\infty]$  
if and only if $A$ is the tensor product of two vectors
in $K_1$.

\bigskip

\noindent REFERENCES
\bigskip

\item{[GVL]} Gene H.~Golub and Charles F.~Van Loan, Matrix
Computations, Johns Hopkins U.~Press, 1983.

\item{[H]} N. J. Higham, Accuracy and Stability of Numerical Algorithms, SIAM
1996.

\item{[HJ]} R.~ A. ~Horn and Charles ~R. Johnson, Matrix Analysis, Cambridge,
1985.

\item{[HLP]} G.~H.~Hardy, J.~E.~Littlewood, and G.~Polya, Inequalities,
Cambridge, 1952.

\item{[SS]} H. Schneider and G.~A.~Strang, Comparison theorems for supremum
\break norms, Numer. Math. 4:15--20, 1962.

\end